\documentclass{AIMS}
\usepackage{amsmath}
  \usepackage{paralist}
  \usepackage{graphics} %% add this and next lines if pictures should be in esp format
  \usepackage{epsfig} %For pictures: screened artwork should be set up with an 85 or 100 line screen
\usepackage{amsthm}
\usepackage{verbatim}
 \usepackage[colorlinks=true]{hyperref}
\hypersetup{urlcolor=blue, citecolor=red}

\theoremstyle{plain}
\newtheorem{thm}{Theorem}[section]
\newtheorem{lem}[thm]{Lemma}
\newtheorem{cor}[thm]{Corollary}
\newtheorem{prop}[thm]{Proposition}

\title[Periods Lengths of PCFGs]
      {On the Period Lengths of the Parallel Chip-Firing Game}

\author[Tian-Yi Jiang]{}


 \email{jiangty@mit.edu}
\thanks{Research supported by the Center for Excellence in Education.}

\begin{document}
\vspace{-2cm}
\maketitle
% Enter the first author's name and address:
\centerline{\scshape Tian-Yi Jiang}
\medskip
{\footnotesize
% please put the address of the first author
 \centerline{}
   \centerline{}
   \centerline{}
} % Do not forget to end the {\footnotesize by the sign }

%\medskip

%\centerline{\scshape First-name2 last-name2 and First-name3
%last-name3}
%\medskip
%{\footnotesize
 % please put the address of the second  and third author
% \centerline{ First line of the address of the second author}
%   \centerline{Other lines}
%   \centerline{Springfield, MO 65810, USA}
%}

\bigskip

% The name of the associate editor will be entered by an editorial staff
% \centerline{(Communicated by the associate editor name)}

%The abstract of your paper
\begin{abstract}
The parallel chip-firing game is a periodic automaton on graphs in which vertices ``fire'' chips to their neighbors.  In 1989, Bitar conjectured that the period of a parallel chip-firing game with $n$ vertices is at most $n$. Though this conjecture was disproven in 1994 by Kiwi et.\ al., it has been proven for particular classes of graphs, specifically trees (Bitar and Goles, 1992) and the complete graph $K_n$ (Levine, 2008). We prove Bitar's conjecture for complete bipartite graphs and characterize completely all possible periods for positions of the parallel chip-firing game on such graphs. 
Furthermore, we extend our construction of all possible periods for games on the bipartite graph to games on complete $c$-partite graphs, $c>2$, and prove some pertinent lemmas about games on general simple connected graphs.
\end{abstract}
%The title of your section 1
\section{Introduction}

\subsection{Definitions}

The parallel chip-firing game or candy-passing game is a periodic automaton on graphs in which vertices, each of which contains some nonnegative number of chips, ``fire'' exactly one chip to each of their neighbors if possible. Formally, let $G$ be an undirected graph with vertex set $V(G)$ and edge set $E(G)$. Define the \emph{parallel chip-firing game} on $G$ to be an automaton governed by the following rules:

\begin{itemize}
\item At the beginning of the game, $\sigma(v)$ chips are placed on each vertex $v$ in $G$, where $\sigma(v)$ is a nonnegative integer. Let a \emph{position} of the parallel chip-firing game, denoted by $\sigma$, be the ordered pair $(G,\{\sigma(v), v \in G\})$ containing the graph and the number of chips on each vertex of the graph.
\item At each \emph{move} or \emph{step} of the game, if a vertex $v$ has at least as many chips as it has neighbors, it will give (\emph{fire}) exactly one chip to each neighbor. Such a vertex is referred to as \emph{firing}; otherwise, it is \emph{non-firing}. All vertices fire simultaneously (in parallel). %the number of chips $\sigma(v)$ on a vertex $v$ is at least the degree of $v$, one chip is fired down each edge incident to $v$, landing on each neighboring vertex of $v$. Each of these firings occur at the same time.
\end{itemize}

We employ the notation of Levine \cite{levine}. Let $U$ denote the step operator; that is, $U\sigma$ is the position resulting after one step is performed on $\sigma$. Let $U^0\sigma=\sigma$, and $U^m\sigma = U U^{m-1}\sigma$. We refer to $U^m\sigma$ as the position occurring \textit{after $m$ steps}. %, and say a vertex $v$ fires \textit{at step $m$} if $U^m\sigma(v) \ge \textrm{deg}(v)$. 
For simplicity, we limit our discussion to connected graphs.

As the number of chips and number of vertices are both finite, there are a finite number of positions in this game. Additionally, since each position completely determines the next position, it follows that for each initial position $\sigma$, there exist some positive integers $p$ such that for large enough $t$, $U^t\sigma = U^{t+p}\sigma$. We refer to the minimal such $p$ as the \emph{period} $p(\sigma)$ of $\sigma$, and we refer to the set $\{U^t\sigma, U^{t+1}\sigma, \ldots, U^{t+p-1}\sigma\}$ as one period of $\sigma$. Also, we call the minimal such $t$ the \emph{transient length} $t_0$ of $\sigma$.

\subsection{Notation}
For easy reference, we include a table of definitions here. We will
focus on the bipartite complete graph for most of the paper, where the
two parts of the graph will be denoted $L$ and $R$. If a notation refers
to $L$, there will be an analogous notation for $R$ as well.

\noindent
\begin{tabular}{lp{3.25in}}
Notation                 & Description \\ \hline
$\sigma$ & the initial position of the game \\
$\sigma(v)$ & chips placed on $v$ at beginning of game  \\
$U^m\sigma$ & position of game after $m$ steps \\
$p(\sigma)$ & period of game starting at $\sigma$ \\
$u_t(\sigma, v)$ & the number of times $v$ fires in the first $t$ steps
\\
$F_v(t)$ & the indicator function of whether a vertex $v$ fires at step
$t$ \\
$F_L(\sigma)$ & the number of vertices that fire in $L$ at position
$\sigma$\\
$d_t(v,m)$ & the total number of times $v$ fires starting from, and
including, $U^m\sigma$, in $t$ steps. $m = 0$ is suppressed. \\
$\alpha_t(L, m)$ & the total number of times vertices in $L$ fire
starting from, and including, $U^m\sigma$, in $t$ steps. $m = 0$ is
suppressed. \\
\end{tabular}

\subsection{Previous Work}

The parallel chip-firing game was introduced by Bitar and Goles \cite{bitar} in 1992 as a special case of the general chip-firing game posited by Bj\"orner, Lov\'asz, and Shor \cite{bjorner} in 1991. They \cite{bitar} showed that the period of any position on a tree graph is 1 or 2. In 2008, Kominers and Kominers \cite{kom1,kom2} further showed that all connected graphs satisfying $\displaystyle\sum_{v \in G} \sigma(v) \ge 4|E(G)| - |V(G)|$ %(where $|E(G)|$ and $|V(G)|$ denote the number of edges and vertices in $G$) 
have period 1; they further established a polynomial bound for the transient length of positions on such graphs. Their result \cite{kom1} that the set of all ``abundant'' vertices $v_i$ with $\sigma(v_i) \ge 2\ \textrm{deg}(v_i)$ stabilizes is particularly useful in simplifying the game.

It was conjectured by Bitar \cite{bitar2} that $p(\sigma) \le |V(G)|$ for all games on all graphs $G$. 
However, Kiwi et.\ al.\ \cite{kiwi} constructed a graph on which there existed a position whose period was at least $\exp(\Omega(\sqrt{|V(G)|\log|V(G)|}))$, disproving the conjecture. Still, it is thought that excluding particular graphs constructed to force long periods, most graphs still have periods that are at most $|V(G)|$. In 2008, Levine \cite{levine} proved this for the complete graph $K_n$. 
%Here, we note some lemmas on chip-firing for general graphs, and . 

\subsection{A Broader Perspective}

The parallel chip-firing game is a special case of the more general chip-firing game, in which at each step, a vertex is chosen to fire. The general chip-firing game, in turn, is an example of an \emph{abelian sandpile} \cite{bjorner}, and has been shown to have deep connections in number theory, algebra, and combinatorics, ranging from elliptic curves \cite{musiker} to the critical group of a graph \cite{biggs} to the Tutte polynomial \cite{lopez}. Bitar and Goles \cite{bitar} observed that the parallel chip-firing game has ``nontrivial computing capabilities,'' being able to simulate the AND, NOT, and OR gates of a classical computer; later, Goles and Margenstern \cite{goles2} showed that it can simulate any two-register machine, and therefore solve any theoretically solvable computational problem. Finally, the parallel chip-firing game can be used to simulate a pile of particles that falls whenever there are too many particles stacked at any point; this important problem in statistical physics is often referred to as the \emph{deterministic fixed-energy sandpile} \cite{casartelli,bagnoli}. The fixed-energy sandpile, in turn, is a subset of the more general study of the so-called \emph{spatially extended dynamical systems}, which occur frequently in the physical sciences and even economics \cite{bak}. Such systems demonstrate the phenomenon of \emph{self-organized criticality}, tending towards a ``critical state'' in which slight perturbations in initial position cause large, avalanche-like disturbances. Self-organized critical models such as the abelian sandpile tend to display properties of real-life systems, such as $1/f$ noise, fractal patterns, and power law distribution \cite{bak, bak2}. Finally, the parallel chip-firing game is an example of a cellular automaton, the study of which have implications from biology to social science.

\subsection{Our Results}

In Section 2, we establish some lemmas about parallel chip-firing games on general simple connected graphs. We bound the number of chips on any single vertex in games with nontrivial period, define the notion of a complement position $\sigma_c$ of $\sigma$ and show that it has the same behavior as $\sigma$, and find a necessary and sufficient condition for a period to occur. Then, in Section 3, we find, with proof, every possible period for the complete bipartite graph $K_{a,b}$. We do so by first showing the only possible periods are of length $k$ or $2k$ for $k \le \min(a,b)$, and then constructing  games with such periods, proving our main result. Finally, in Section 4, we construct positions on the complete $c$-partite graph $K_{a_1, a_2, \ldots, a_c}$ with period $p$ for all $1 \le p \le c \cdot \min(a_1, a_2, \ldots, a_c)$.
%The title of your section 2

\section{Parallel Chip-Firing on Simple Connected Graphs}

Consider a simple connected graph $G$. For each vertex $v$ in $G$, let $\Phi_{\sigma}(v)$ denote the number of firing neighbors $w$ of $v$; that is, the number of vertices $w$ neighboring $v$ satisfying $\sigma(w) \ge \textrm{deg}(w)$. A step of the parallel chip-firing game on $G$ is then defined as follows:

\begin{equation}
\label{eq:firing}
U\sigma(v) = 
\begin{cases}
\sigma(v) + \Phi_{\sigma}(v), & \sigma(v) \le \textrm{deg}(v) - 1\\
\sigma(v) + \Phi_{\sigma}(v) - \textrm{deg}(v), & \sigma(v) \ge \textrm{deg}(v).
\end{cases}
\end{equation}

Define a \emph{terminating} position to be a position in which no vertices fire after finitely many moves. We begin our investigation by proving some lemmas limiting the number of chips on each vertex in a game with nontrivial period (period greater than 1). %Our first lemma, due to Kominers and Kominers \cite{kom1}, bounds the number of chips on any single vertex in a game of nontrivial period.

%%%%%%%%%%% lemma lem:2n-1. lem 2.1 %%%%%%%%%%%%%
\begin{lem} \label{lem:2n-1} For sufficiently large $t$, $U^t\sigma(v) \le 2\ \textup{deg}(v)-1$ for all $v \in G$ in all games with nontrivial period on a connected graph $G$. \end{lem}%unless all vertices fire at each step.\end{lem}

\begin{proof} %Notice that if a $\sigma(v) \ge \textrm{deg}(v)$, then $U\sigma(v) \le \textrm{deg}(v)$ because $v$ fires $\textrm{deg}(v)$ chips but receives at most $\textrm{deg}(v)$ chips. Assume for contradiction that there exists some $v \in G$ for which $\sigma(v) < 2\ \textrm{deg}(v)$ but $U\sigma(v) \ge 2\ \textrm{deg}(v)$. Because $v$ may receive at most $\textrm{deg}(v)$ chips, $\sigma(v) > \textrm{deg}(v)$; but then $v$ will fire, meaning $U\sigma(v) \le \sigma(v) < 2\ \textrm{deg}(v)$, a contradiction. Hence, for sufficiently large $t$, we may assume that if a vertex $v \in G$ satisfies $U^t\sigma(v) \ge 2\ \textrm{deg}(v)$, then $U^{t+1}\sigma(v) = U^t\sigma(v)$, because the number of chips on $v$ is nonincreasing. 
Kominers and Kominers \cite{kom1} showed that if a vertex $v \in G$ satisfies $\sigma(v) \le 2\ \textrm{deg}(v)-1$, then $U\sigma(v) \le 2\ \textrm{deg}(v)-1$. % no vertex with fewer than $2\ \textrm{deg}(v)$ chips after any one step may contain at least $2\ \textrm{deg}(v)$ chips after the next step. 
They then showed that if, after sufficiently many steps $t$, there still exists a vertex $v$ with $U^t\sigma(v) \ge 2\ \textrm{deg}(v)$ , then all vertices must be firing from that step onward. Since the period of a position is 1 if and only if either all or no vertices in $G$ are firing \cite{bitar}, $U^t\sigma(v) \le 2\ \textrm{deg}(v)-1$ is true for any game on $G$ with nontrivial period and sufficiently large $t$.%Hence, the lemma is true.  
\end{proof} %the set of vertices with $\sigma(v)\ge 2\ \textrm{deg}(v)$ eventually , and if one exists in the period, all vertices must be firing. 

We further bound the number of chips on each vertex by generalizing a result of Levine \cite{levine}:

%%%%%%%%%%% lemma lem:confined. lem 2.2 %%%%%%%%%%%%%
\begin{lem} \label{lem:confined}
Consider a vertex $v$ in position $\sigma$ such that $\sigma(v)\le 2\ \textrm{deg}(v)-1$. Then $$\Phi_{\sigma}(v) \le U\sigma(v) \le \Phi_{\sigma}(v)+\textup{deg}(v)-1.$$  \end{lem}

\begin{proof} 
Either $\sigma(v) < \textrm{deg}(v)$ or not. We consider the cases individually.

If $0\le\sigma(v)\le \textrm{deg}(v)-1$, then $U\sigma(v) = \sigma(v)+\Phi_{\sigma}(v)$. So 
$$\Phi_{\sigma}(v) \le U\sigma(v) \le \Phi_{\sigma}(v) + \textrm{deg}(v)-1.$$

If instead $\textrm{deg}(v) \le \sigma(v) \le 2\ \textrm{deg}(v)-1$, then $U\sigma(v) = \sigma(v)+\Phi(v)-\textrm{deg}(v)$. Hence 
$$\Phi(v) = \textrm{deg}(v)+\Phi_{\sigma}(v)-\textrm{deg}(v) \le \sigma(v)+\Phi_{\sigma}(v)-\textrm{deg}(v)$$
$$= U\sigma(v) \le 2\ \textrm{deg}(v)-1 + \Phi_{\sigma}(v) - \textrm{deg}(v) = \textrm{deg}(v)-1+\Phi_{\sigma}(v).\mbox{\qedhere} $$ \end{proof}

If a vertex $v$ satisfies $\Phi_{\sigma}(v) \le \sigma(v) \le \Phi_{\sigma}(v)+\textrm{deg}(v)-1$, we call it \emph{confined}. Furthermore, call a position confined if all vertices in the position are confined.  Note that for confined $v$, $\sigma(v) \le \Phi_{\sigma}(v)+\textrm{deg}(v)-1 \le 2\ \textup{deg}(v)-1.$ Lemmas ~\ref{lem:2n-1} and~\ref{lem:confined} imply that if $p(\sigma)>1$, then $U^t\sigma$ is confined if $t \ge t_0$, where $t_0$ is the transient length of $\sigma$; that is, once the game reaches a position which repeats periodically, all subsequent positions are confined. We generally limit our discussion to confined positions to exclude positions with trivial periods.

Next, we define 
\begin{equation} \label{eq:fv}
F_v(t) = 
\begin{cases}
1, & U^t\sigma(v) \ge \textrm{deg}(v)\\
0, & U^t\sigma(v) \le \textrm{deg}(v)-1
\end{cases} 
\end{equation}
to be the indicator function of whether a vertex $v$ fires at step $t$. We prove a lemma about positions that are equivalent, or have the same behavior, when acted upon by the step operator $U$.

%%%%%%%%%%% lemma lem:complement %%%%%%%%%%%%%
\begin{lem} \label{lem:complement} Let the \emph{complement} $\sigma_c$ of a confined position $\sigma$ be the position that results after replacing the $\sigma(v)$ chips on each vertex $v\in G$ with $2\ \textrm{deg}(v)-1-\sigma(v)$ chips. Then $U(\sigma_c) = (U\sigma)_c$. \end{lem}

\begin{proof} We begin by noticing that since $\sigma$ is confined, each vertex $v$ has at most $2\ \textrm{deg}(v)-1$ chips, so each vertex in $\sigma_c$ has a nonnegative number of chips.

Observe that a vertex $v$ fires in $\sigma_c$ exactly when it did not fire in $\sigma$. Hence, $U\sigma(v) = \sigma(v) + \Phi_{\sigma}(v) - F_v(0)\textrm{deg}(v)$, and all but $\Phi_{\sigma}(v)$ neighbors will fire in $\sigma_c(v)$. So
%$$U(\sigma_c)(v) = (2\ \textrm{deg}(v)-1 - \sigma(v)) + (\textrm{deg}(v) - \Phi_{\sigma}(v)) = (2\ \textrm{deg}(v)-1)-(\Phi_{\sigma}(v)+\sigma(v)-\textrm{deg}(v)) = (U\sigma)_c(v)$$
\begin{align*}
U(\sigma_c(v))&=(2\ \textrm{deg}(v)-1-\sigma(v))+(\textrm{deg}(v)-\Phi_{\sigma}(v))-((1-F_v(0))\textrm{deg}(v))\\
&=(2\ \textrm{deg}(v)-1)-(\Phi_{\sigma}(v)+\sigma(v)-F_v(0)\textrm{deg}(v)) = (U\sigma)_c(v). \mbox{\qedhere}
\end{align*} 
\end{proof}

This lemma means we may treat $\sigma$ and $\sigma_c$ as equivalent positions, as at any point during their firing, we may transform one into the other. This implies the following corollary:

\begin{cor} For all positions $\sigma$ on $G$, $p(\sigma) = p(\sigma_c)$.\end{cor}

 Next, we prove a proposition that characterizes a period of the game on any connected graph $G$.
For each position $\sigma$ and vertex $v\in G$, let
 $$u_t(\sigma,v) = \Bigg|\; \{s \; | \; 0\le s < t, U^s\sigma(v) \ge \textrm{deg}(v)\} \;\Bigg|.$$
be the number of times $v$ fires in the first $t$ steps.

%%%%%%%%%%% proposition prop: alleq %%%%%%%%%%%%%
\begin{prop} \label{prop:alleq} The position $\sigma$ on $G$ satisfies $U^t\sigma = \sigma$ if and only if each vertex has fired the same number of times within those $t$ steps; that is, iff for all vertices $v,w\in G$,
\begin{equation} 
\label{eq:alleq}
u_t(\sigma,v) = u_t(\sigma,w) = k \ge 0.
\end{equation}
\end{prop}
%The only if portion of the following proof is due to Brian Hamrick.
\begin{proof} 
If equation 
\eqref{eq:alleq} holds, then by equation \eqref{eq:firing}, $U^t\sigma(v) = \sigma(v) + k \cdot \textrm{deg}(v) - k \cdot \textrm{deg}(v) = \sigma(v)$ for all $v$, so $U^t\sigma = \sigma$. Conversely, if $U^t\sigma = \sigma$, consider the vertex $v'$ such that $u_t(\sigma,v')=k'$ is maximal. %$ = \max_{v\in G} u_t(\sigma,v)=k'$. 
Then, since $u_t(\sigma,w) \le k'$ for all vertices $w$ neighboring $v$, 
$$U^t\sigma(v) = \sigma(v) + \sum_{w}u_t(\sigma,w) - k'\  \textrm{deg}(v) \le \sigma(v) + k'\ \textrm{deg}(v) - k'\ \textrm{deg}(v).$$
 But as $U^t\sigma(v) = \sigma(v)$, we see that $u_t(\sigma,w) = k'$ must hold for all $w$ neighboring $v$. Since the graph is connected, we continue inductively through the entire graph to obtain equation \eqref{eq:alleq}.
\end{proof}

%SECTION 3 SECTION 3 SECTION 3!
%SECTION 3 SECTION 3 SECTION 3!	
\section{Period Length of Games on $K_{a,b}$}

Recall that a complete bipartite graph $G=K_{a,b}$ may be partitioned into two subsets of vertices, $L$ and $R$, such that no edges exist among vertices in the same set, but every vertex in $L$ is connected to every vertex in $R$. We refer to the sets $L,R$ as the \emph{sides} of $G$. Define $a = |L|$ and $b = |R|$. As stated above, Bitar and Goles \cite{bitar} showed that if no vertices or all vertices are firing, the period is 1. We consider only games whose period is greater than 1; that is, at least one vertex is firing every turn, and not all vertices fire every turn. 

Let $F_L(\sigma)$ and $F_R(\sigma)$ denote the number of vertices in $L$ and $R$, respectively, that fire in $\sigma$.
Then, $\Phi_{\sigma}(v) = F_R(\sigma)$ if $v \in L$, and $\Phi_{\sigma}(v) = F_L(\sigma)$ if $v\in R$. Notice that $F_R(\sigma)$ is the number of vertices in $R$ with at least $a$ chips, and $F_L(\sigma)$ is the number of vertices in $L$ with at least $b$ chips. %on $K_{n,n}$, the parallel chip-firing game may be defined as follows:
 Let $\alpha_t(L,m) = \sum_{v\in L} (u_{m+t}(\sigma,v)-u_m(\sigma,v))$ be the number of times any of the vertices in $L$ have fired in the first $t$ steps starting from, and including, $U^m\sigma$, and define $\alpha_t(R,m)$ similarly. Define $\alpha_t(L) = \alpha_t(L,0)$ and $\alpha_t(R) = \alpha_t(R,0)$. 
 
Without loss of generality, we prove facts about the vertices in $L$, which also hold for vertices in $R$. In the first $t$ steps, a vertex $v$ in $L$ fires a total of $b u_t(\sigma,v)$ chips and receives $\alpha_t(R)$ chips. Hence,

\begin{equation}
\label{eq:firediff}
U^t\sigma(v) - \sigma(v) = \alpha_t(R) - b u_t(\sigma,v).
\end{equation}

Next, we prove a lemma that bounds the number of times a vertex has fired once the position is confined.

%%%%%%%%%%% lemma lem:diff1 %%%%%%%%%%%%%
\begin{lem} \label{lem:diff1} Let $v,w \in L$. If $\sigma$ is confined, and $\sigma(v) \le \sigma(w)$, then for all $t \ge 0$,
$$u_t(\sigma,v) \le u_t(\sigma,w) \le u_t(\sigma,v) + 1.$$ \end{lem}

\begin{proof}
We prove this by induction on $t$. The base case, $t=1$, is straightforward: vertices $v$ and $w$ have each fired either 0 or 1 times. If $v$ fires after step $0$, then $\sigma(w) \ge \sigma(v) \ge \textrm{deg}(v) = \textrm{deg}(w)$ chips, and $w$ also fires. Now, assume $u_t(\sigma,v) \le u_t(\sigma,w) \le u_t(\sigma,v) + 1.$

If $u_t(\sigma,w) = u_t(\sigma,v)$, then by equation \eqref{eq:firediff},
$$U^t\sigma(w) - \sigma(w) - (U^t\sigma(v) - \sigma(v)) = \alpha_t(R) - b u_t(\sigma,w) - (\alpha_t(R) - b u_t(\sigma, v))$$
$$\implies U^t\sigma(w) - U^t\sigma(v) = \sigma(w) - \sigma(v) \ge 0.$$

Thus, if $v$ is ready to fire after step $t$, then $w$ must be ready to fire also. It follows that 
\begin{equation}
\label{eq:firediff2}
u_{t+1}(\sigma,v) \le u_{t+1}(\sigma,w) \le u_{t+1}(\sigma,v) + 1.
\end{equation}

Otherwise, $u_t(\sigma,w) = u_t(\sigma,v)+1$. Then, since $U^t\sigma$ is confined from Lemma~\ref{lem:confined}, by equation \eqref{eq:firediff}, 
$$U^t\sigma(v) - \sigma(v) - (U^t\sigma(w) - \sigma(w)) = \alpha_t(R) - b u_t(\sigma,v) - (\alpha_t(R) - b u_t(\sigma, w))$$
$$\implies U^t\sigma(v) - U^t\sigma(w) = b +\sigma(v) - \sigma(w) \ge 0$$
by Lemma~\ref{lem:confined}, since the degrees of both $v$ and $w$ are $b$.

So, if $w$ is ready to fire after step $t$, so is $v$, and equation \eqref{eq:firediff2} again holds.
\end{proof}

%\begin{lem}  If $\sigma$ is confined, $p(\sigma) > 1$, and $F_t(L) = 0$ or $n$, then $F_{t+1}(R) = F_t(L)$. \end{lem}

%\begin{proof} From (1.3), it suffices to prove this statement if $F_t(L) = 0$; taking the complement of $\sigma$ for this case yields the case when $F_t(L) = n$. Consider $v \in R$. If $v$ fires at step $t$, then $U^{t+1}\sigma(v) = U^t\sigma(v) + F_t(L) - n = U^t\sigma(v) - n \le 2n-1 - n = n-1$ from %ef{lem:confined}. Hence $v$ does not fire at step $t+1$. If $v$ does not fire at step $t$, it gains no chips because $F_t(L) = 0$, so $v$ does not fire at step $t+1$. Hence, $F_{t+1}(R) = 0 = F_t(L)$ as desired. \end{proof}

%This lemma states that if every or no vertex on one side fires at step $t$, then at step $t+1$, every or no vertex, respectively, fires on the other side.

%\begin{cor} \label{cor:bothfire} If $\sigma$ is confined, $n > F_t(L) > 0$ and $n > F_t(R) > 0$, then $n > F_{t+1}(L) > 0$ and $n > F_{t+1}(R) > 0$. \end{cor}

%This corollary divides the set of positions $\sigma$ into those that have some, but not all, vertices on both sides firing, and those that have either no or all vertices on one side firing at any step.

From the above lemma, we can deduce the following: %The following lemma yields more information on the number of times an individual vertex has fired given the total number of times an entire side has fired.

%%%%%%%%%%% lemma lem:ndiv %%%%%%%%%%%%%
\begin{lem} \label{lem:ndiv} If $\sigma$ is confined and $a|\alpha_t(L)$, then $u_t(\sigma,v) = \alpha_t(L)/a$ for all $v\in L$. \end{lem}

\begin{proof} Let $v'$ be the vertex in $L$ with $\sigma(v')$ minimal. 
By Lemma~\ref{lem:diff1}, for all $v \in L$, 
$$u_t(\sigma,v) \in \{m, m+1\}$$ where $m = u_t(\sigma,v')$. 
If $z$ is the number of vertices $w \in L$ with $u_t(\sigma,w) = m+1$, then  
$$\alpha_t(L) = \sum_{v\in L} u_t(\sigma,v) = (a-z)m + z(m+1) = am+z \equiv z \pmod{a}.$$  
Since $u_t(\sigma,v') = m$, we have $z < a$. Then $z = 0$ because $z \equiv 0 \pmod{a}$; so $u_t(\sigma,v) = m$ for all $v \in L$. Since $\alpha_t(L) = \sum_{v\in L}(u_t(\sigma,v))$, this implies $u_t(\sigma,v) = \alpha_t(L)/a$ for all $v \in L$. \end{proof}

Clearly, if all $a$ vertices in $L$ have fired the same number of times, then $a|\alpha_t(L)$; so we have found a necessary and sufficient condition for all vertices on the same side to fire the same number of times. But by Proposition~\ref{prop:alleq}, a period is completed when all $v\in G$ have fired the same number of times; thus, we desire a relation between the sides that forces every vertex on both sides to fire the same number of times. Our first step is the following lemma. 

%%%%%%%%%%% lemma lem:t1 %%%%%%%%%%%%%
\begin{lem} 
\label{lem:t1} If $\sigma$ is confined, and $\alpha_t(L) = ka$ for some nonnegative integer $k$, then $u_{t+1}(\sigma,v) - u_1(\sigma,v) = k$ for all $v \in R$. \end{lem}
\begin{proof} 
%Let $\sigma_L$ denote the configuration of chips on $L$ in position $\sigma$, and define $\sigma_R$ similarly. 

%Classify each vertex $v$ in $R$ as \emph{firing} if $\sigma(v) \ge a$ and \emph{non-firing} otherwise.

If $v\in R$ is firing, then $U^{t}\sigma(v) = (\sigma(v) - a) + ka - a (u_t(\sigma,v)-1).$ Since $\sigma$ is confined and $v$ is firing, $0\le \sigma(v) - a \le a-1$; and since $U^{t}\sigma$ is confined by Lemma~\ref{lem:confined}, we have $0 \le \sigma(v) - a + ka - a (u_t(\sigma,v)-1) \le 2a-1.$
These two inequalities together imply that, for firing vertices $v$,
$$-a < -(a-1) \le ka - a (u_t(\sigma,v) - 1) = ka - a (u_t(\sigma,v) - u_1(\sigma,v)) \le 2a-1 < 2a.$$  

If $v \in R$ is instead non-firing, then $U^{t}\sigma(v) = \sigma(v) + ka - a u_t(\sigma,v)$ chips. $U^t\sigma(v)$ is confined by Lemma~\ref{lem:confined}, so $0\le\sigma(v) + ka - a u_t(\sigma,v)\le 2a-1$; since $0 \le \sigma(v) \le a-1$ because $v$ is non-firing, we then deduce, similarly as above, that
$$-a < -(a-1) \le ka - a u_t(\sigma,v) = ka - a (u_t(\sigma,v) - u_1(\sigma,v)) \le 2a-1 < 2a$$
for non-firing vertices $v$ as well. 

Therefore, for all $v \in R$, we have that
$$-a < a( k - (u_t(\sigma,v) - u_1(\sigma,v)) ) < 2a \implies -1 <  k - (u_t(\sigma,v) - u_1(\sigma,v))  < 2, $$
so $u_t(\sigma,v)-u_1(\sigma,v) \in \{k, k-1\}$ for all $v \in R$. If $u_t(\sigma,v) = k$, then we can compute $U^{t}\sigma(v) < a = \textrm{deg}(v)$; hence $v$ does not fire after step $t$, and $u_{t+1}(\sigma,v) - u_1(\sigma,v) = u_t(\sigma,v) - u_1(\sigma,v) = k$. If instead $u_t(\sigma,v) = k-1$, then $U^t\sigma(v) \ge a$, so $v$ fires after step $t$, and $u_{t+1}(\sigma,v) - u_1(\sigma,v) = u_t(\sigma,v) - u_1(\sigma,v) + 1 = k$, and we are done.
\end{proof}

Note that applying this lemma to $U^m\sigma$ also means $\alpha_t(L,m) = ka \implies \alpha_t(R,m+1) = kb$.
%Next, define
%\begin{equation}
%F_v(t) = 
%\begin{cases}
%1, & U^t\sigma(v) \ge \textrm{deg}(v)\\
%0, & U^t\sigma(v) \le \textrm{deg}(v)-1,
%\end{cases} 
%\end{equation}

Next, recalling the definition of $F_v(t)$ in equation \eqref{eq:fv}, we define $d_t(v,m) = u_{m+t}(\sigma,v) - u_{m}(\sigma,v) = \sum_{i=m}^{m+t-1}F_v(i)$ for nonnegative integers $m$ and positive integers $t$. Note that by definition, $\sum_{v \in L} d_t(v,m) = \alpha_t(L,m).$ Applying Lemmas~\ref{lem:ndiv} and~\ref{lem:t1} to the position $U^m\sigma$, we find that if
$\alpha_t(L,m) = ka$, then 
\begin{equation} \label{eq:t1}
d_t(v,m) =  d_{t}(w,m+1) = k
\end{equation}
 for all vertices $v \in L$ and all vertices $w \in R$. 

By applying Lemma~\ref{lem:t1} to $L$ (at step 0) and then $R$ (at step 1), we obtain this corollary:
%%%%%%%%%%% corollary cor:t1 %%%%%%%%%%%%%
\begin{cor}\label{cor:t1} If $a|\alpha_t(L,m)$, then $d_t(v,m) = d_t(v,m+2) = k$ for all $v \in L$ and $a|\alpha_t(L,m+2)$. \end{cor}

Now, we give a sufficient condition for a period of a position on $K_{a,b}$ to occur.

%%%%%%%%%%% lemma lem:periodt %%%%%%%%%%%%%
\begin{lem} \label{lem:periodt} Let $k$ be a nonnegative integer. If $\sigma$ is confined, and for some $m \ge 0$ and $t \ge 1$, $\alpha_t(L,m) = ka$ and $F_v(m) = F_v(m+t)$ for all $v \in L$, then $U^t(\sigma) = \sigma$. \end{lem}
\begin{proof}
By Lemma~\ref{lem:ndiv} applied to $U^m(\sigma)$, since $\alpha_t(L,m) = ka$, $d_t(v,m)=k$ for all $v \in L$. 
Since $F_v(m) = F_v(m+t)$ for all $v \in L$,  
$$k=d_t(v,m) = F_v(m) + \displaystyle\sum_{i=m+1}^{m+t-1} F_v(i) = F_v(m+t) + \sum_{i=m+1}^{m+t-1} F_v(i) = d_t(v,m+1)$$
 for all $v \in L$. But by equation~\eqref{eq:t1}, $d_t(v,m) = d_t(w,m+1)=k$ for all $v \in L$, $w \in R$. Hence, $d_t(v,m+1) = d_t(w,m+1)=k$ for all vertices $v,w \in G$, which by Proposition~\ref{prop:alleq} applied to $U^{m+1}\sigma$ implies $U^{m+1+t}\sigma = U^{m+1}\sigma$, or $p(\sigma)|t$. But $t$ is taken to be as small as possible, so $p(\sigma) = t$. 
\end{proof}

We can strengthen this lemma slightly:
%%%%%%%%%%%% corollary cor:periodt %%%%%%%%%%%%
\begin{cor} \label{cor:periodt} Let $k$ be a nonnegative integer. If $\sigma$ is confined, and for some $m \ge 0$ and $t \ge 1$, $\alpha_t(L,m) = ka$ and $F_L(m) = F_L(m+t)$, then $U^t(\sigma) = \sigma$. \end{cor}
\begin{proof} It suffices to show that $F_L(m) = F_L(m+t)$ and $a|\alpha_t(L,m)$ together imply $F_v(m) = F_v(m+t)$ for every $v \in L$. Note that since $a|\alpha_t(L,m)$, then by Lemma~\ref{lem:ndiv}, $d_t(v,m) = k$ for all $v \in L.$ But from steps $t$ to $t+m$, each of these $v$ has received the same number of chips from $R$, $k_0$, while each firing the same number of chips $k$ to $R$. Hence 
\begin{equation}\label{eq:samefirediff}U^{m+t}(\sigma(v)) = U^m(\sigma(v)) + (k_0-k).\end{equation}

Sort the vertices in $L$ by their chip numbers at step $m$, and call these vertices $v_1, v_2, \ldots, v_a$ so that $i \le j \Leftrightarrow U^m(\sigma(v_i)) \le U^m(\sigma(v_j)).$ Hence the set of firing vertices $S_m$ at step $m$ must be the set $\{v_i, i \ge c_0\}$ for some $1 \le c_0 \le a.$ 

It follows from \eqref{eq:samefirediff} that $U^{m+t}(\sigma(v_i)) = U^m(\sigma(v_i)) + (k_0-k) \le U^m(\sigma(v_j)) + (k_0-k) = U^{m+t}(\sigma(v_j)).$ Again, the set of firing vertices $S_{m+t}$ at step $m+t$ must be the set $\{v_i, i \ge c_1\}$ for some $1 \le c_1 \le a.$ But since $|S_m| = |S_{m+t}|$, we conclude that $c_0 = c_1$ and hence $S_m = S_{m+t}.$ 
\end{proof}

With this corollary, we show the following lemma:
%%%%%%%%%% proposition prop:t1t2 %%%%%%%%%%%%5
\begin{prop}\label{prop:t1t2} Let $\sigma$ be a confined position. Let $t = t_1, t = t_1 + t_2$ be the two smallest positive integer solutions to $\alpha_t(L,0) \equiv 0 \pmod{a}.$ Then $p(\sigma) | t_1+t_2$ and $t_1+t_2$ is even. \end{prop}
\begin{proof}
We first show that $\tau = t_1 + t_2$ is even.

Let $g(m)$ be the smallest positive integer with $\alpha_{g(m)}(L, m) \equiv 0 \pmod{a}.$ Without loss of generality, choose $t_1 = \min(g(m)), m \ge 0.$ Such a $t_1$ exists since $g(m) \ge 0.$ If $\alpha_{t_1}(L, m_1) \equiv 0 \pmod{a}$, we can shift our reference position and call $U^{m_1}(\sigma)$ our new $\sigma.$

We want to show that $g(0) = t_1, g(t_1) = t_2$.

Say $t_1$ is even. Assume for contradiction that $t_2 > t_1$. Then by repeatedly applying Corollary~\ref{cor:t1}, $\alpha_{t_1}(L, t_1) = 0$. But this contradicts $g(t_1) = t_2$. Hence $t_1 = t_2$. 

Now say $t_2$ is even. Assume for contradiction that $t_2 > t_1$. Then by repeatedly applying Corollary~\ref{cor:t1}, $\alpha_{t_1}(L, t_2) \equiv 0 \pmod{a}$. But then $\alpha_{t_2}(L, t_1) - \alpha_{t_1}(L, t_2) = \sum_{i=t_1}^{t_1+t_2-1} F_L(t) - \sum_{i=t_2}^{t_1+t_2-1} F_L(t) = \sum_{i=t_1}^{t_2-1} F_L(t) = \alpha_{t_2-t_1}(L, t_1) \equiv 0 \pmod{a}.$ But since $t_2 - t_1 < t_2$, this again contradicts $g(t_1) = t_2$. Hence $t_1 = t_2$. 

So if either of $t_1, t_2$ are odd, the other is odd as well. Hence $\tau = t_1+t_2$ is even.

Now by Corollary~\ref{cor:t1}, $\alpha_{\tau}(L,m) = \alpha_{\tau}(L,m+2) = k$ for all $k$ even. Then
$$\alpha_{\tau}(L,m) = \alpha_{\tau}(L,m+2) \Leftrightarrow \sum_{i=m}^{m+\tau-1}F_L(i) = \sum_{i=m+2}^{m+\tau+1}F_L(i)$$ 
\begin{equation}\label{eq:twoequal} \Leftrightarrow F_L(m) + F_L(m+1) = F_L(m+\tau) + F_L(m+\tau+1).\end{equation}
Similarly, $F_R(m+1) + F_R(m+2) = F_R(m+\tau+1) + F_R(m+\tau+2)$ by Lemma~\ref{lem:t1} and Corollary~\ref{cor:t1}.

By Corollary~\ref{cor:periodt}, it suffices to show that $F_L(m+ n\tau) = F_L(m+ (n+1)\tau)$ or that $F_R(m+ n\tau + 1) = F_R(m+ t_1 + (n+1)\tau + 1)$ for some nonnegative integer $n$. Since $F_L(t) \le a$, we can find a maximal $F_L(m+ n\tau)$; call this $F_L(m+n_0\tau).$ 

Now assume for contradiction that $ F_L(m + (n_0+1)\tau + k) \neq F_L(m + n_0\tau + k)$ for all even $k \ge 0$. 
We will show that for all even $k \ge 0$, 
\begin{equation}\label{eq:descent}F_L(m+(n_0+1)\tau+k) < F_L(m+n_0\tau+k).\end{equation} 

Since $F_L(m+n_0\tau)$ is maximal,  $F_L(m+(n_0+1)\tau) < F_L(m+n_0\tau)$. Then $\alpha_{\tau}(L, m+(n_0+1)\tau) < \alpha_{\tau}(L, m+n_0\tau)$. Now consider a vertex $w \in R$. From steps $m+n_0\tau+1$ to $m+(n_0+1)\tau$, $w$ has fired $k$ times (since $d_{\tau}(w, m+n_0\tau+1) = \frac{1}{a} \cdot \alpha_{\tau}(L, m+n_0\tau) = k$), but has received a chip only $\alpha_{\tau}(L, m+(n_0+1)\tau) < \alpha_{\tau}(L, m+n_0\tau) = ka$ times. 

So for each $w \in R$, $U^{m+(n_0+1)\tau+1}(\sigma(w)) - U^{m+n_0\tau+1}(\sigma(w)) = \alpha_{\tau}(L, m+(n_0+1)\tau) - ka < 0$, so $U^{m+(n_0+1)\tau+1}(\sigma(w)) < U^{m+n_0\tau+1}(\sigma(w))$ for each $w \in R$. Hence $F_R(m+(n_0+1)\tau) < F_R(m+n_0\tau+1).$

But applying the same logic with $L$ and $R$ reversed, it follows that $F_L(m+(n_0+1)\tau+2) < F_R(m+n_0\tau+2).$ So by induction, our claim \eqref{eq:descent} holds. But since $\tau$ is even, we get an infinitely decreasing sequence of integers $F_L(m+n_0\tau) > F_L(m+(n_0+1)\tau) > F_L(m+(n_0+2)\tau) > \ldots$, a contradiction. Hence for some even $k \ge 0$, $F_L(m + n_0\tau + k) = F_L(m + (n_0+1)\tau + k)$. But by Corollary~\ref{cor:t1}, $\alpha_{\tau}(L, m+n_0\tau+k) \equiv 0 \pmod{a}$, so by Corollary~\ref{cor:periodt}, it follows that $U^{m + (n_0+1)\tau + k}\sigma = U^{m+n_0\tau+k}\sigma$. Hence $\sigma$ eventually repeats every $\tau = t_1+t_2$ steps, so $p(\sigma) | t_1+t_2$ as desired.
\end{proof}

%%%%%%%%%%% corollary cor:lea %%%%%%%%%%%%%
\begin{cor}\label{cor:lea} Let $\sigma$ be a position on $K_{a,b}$. If $p(\sigma)$ is odd, $p(\sigma) \le \min(a,b)$; and if $p(\sigma)$ is even, $p(\sigma) \le 2\min(a,b)$.\end{cor}
\begin{proof} Without loss of generality, let $a \le b$. If $p(\sigma)=1$, $p(\sigma) \le a$. Otherwise, $p(\sigma)>1$. If $t_0$ is the transient length of $\sigma$, then since $p(\sigma) = p(U^{t_0}\sigma)$, we may replace $\sigma$ by $U^{t_0}\sigma$ and assume $\sigma$ is confined by Lemma~\ref{lem:confined}. 

By the Pigeonhole Principle, among every $2a+1$ steps, there are three steps $m, m + t_1, m + t_1 + t_2$ with $\alpha_m(L) \equiv \alpha_{m + t_1}(L) \equiv \alpha_{m + t_1 + t_2}(L) \pmod{a}$. But this means $\alpha_{t_1}(L, t) \equiv \alpha_{t_2+t_1}(L, t) \equiv 0 \pmod{a}$. Also, since $\sigma$ is confined, every position afterward is confined. Hence we can apply Proposition~\ref{prop:t1t2} and conclude that $p(\sigma) | t_1+t_2$. But we have $t_1+t_2 \le 2a+1 - 1 = 2a$, so $p(\sigma)$ divides some even number $t_1+t_2 \le 2a.$ We conclude that if $p(\sigma)$ is odd, then $p \le a$, and if $p(\sigma)$ is even, then $p \le 2a$, as desired.

%By the Pigeonhole Principle, there must exist steps $t_1 \ge t_2 \ge t_1 - a$ with $\alpha_{t_1}(L,0) \equiv \alpha_{t_2}(L,0) \pmod{a}$. But then $a|\alpha_{t_1 - t_2}(L,t_2)$, so by Proposition~\ref{prop:t2t} applied to $U^{t_2}\sigma$, $p(\sigma) = k$ or $2k$, where $k = t_1-t_2 \le a.$ Hence, if $p(\sigma)$ is odd, $p(\sigma) \le a$, and if $p(\sigma)$ is even, $p(\sigma) \le 2a$. 
 \end{proof}

Finally, we characterize all possible periods for $\sigma$.

%%%%%%%%%%% proposition prop:existpos %%%%%%%%%%%%%
\begin{prop} \label{prop:existpos} There exist positions $\sigma$ on $G=K_{a,b}$ with period $k$ and $2k$ for all $1 \le k \le \min(a,b)$.\end{prop}
\begin{proof}
 Without loss of generality, let $a \le b$.

Let $L_1, L_2, \ldots, L_a$ be the vertices in $L$, and $R_1, R_2, \ldots, R_b$ be the vertices in $R$. Let $k$ be a positive integer such that $2 \le k \le a$. We represent each position $\sigma$ on $G$ by two vectors 
$$L(t) = \Bigg(U^t\sigma(L_1), U^t\sigma(L_2), U^t\sigma(L_3), \ldots, U^t\sigma(L_a)\Bigg),$$
$$R(t) = \Bigg(U^t\sigma(R_1), U^t\sigma(R_2), U^t\sigma(R_3), \ldots, U^t\sigma(R_b)\Bigg).$$

Consider the following position $\sigma_k$, which we claim has period $k$: % where $\sigma(L_i) = \sigma(R_i) = i$
%for $1 \le i \le k-1$ and $\sigma(L_j) = b, \sigma(R_{j'}) = a$
%for all k \le j \le a, k \le j' \le b$.
\[ L(0) = \Bigg(1, 2, \ldots, k-2, k-1, \underbrace{b, b, \ldots, b}_{(a-k+1)\ \textrm{times}} \Bigg), \]
\[ R(0) = \Bigg(1, 2, \ldots, k-2, k-1, \underbrace{a, a, \ldots, a}_{(b-k+1)\ \textrm{times}} \Bigg). \]
%%%

 %The above position is represented as 
$L_{k}, L_{k+1}, \ldots L_a$ and $R_{k}, R_{k+1}, \ldots R_b$ fire, so $U\sigma_k$ is represented by
\[ L(1) = \Bigg(b-k+2, b-k+3, \ldots, b-1, b, \underbrace{b-k+1, b-k+1, \ldots, b-k+1}_{(a-k+1)\ \textrm{times}}\Bigg), \]
\[ R(1) = \Bigg(a-k+2, a-k+3, \ldots, a-1, a, \underbrace{a-k+1, a-k+1, \ldots, a-k+1}_{(b-k+1)\ \textrm{times}}\Bigg). \]
%%%
We can see that the vertices $L_i, R_j$ with $U^t\sigma_k(L_i) = b, U^t\sigma_k(R_i) = a$ satisfy $i = k-t$ for $t = 1, 2, \ldots, k-1$. So, $U^t\sigma_k = \sigma_k$ follows upon applying Proposition~\ref{prop:alleq}, noting that after $k$ steps, each vertex has fired exactly once. Hence, $\sigma_k$ has period $k$.

Next, consider the following position $\sigma_{2k}$, which we claim has period $2k$: %where $\sigma(L_i) = i-1, \sigma(R_i) = i$ for $1 \le i \le k-1$
%and $\sigma(L_j) = k-1, \sigma(R_{j'}) = a$ for $k \le j \le a, k \le j' \le b$.
\[ L(0) = \Bigg(0, 1, \ldots, k-3, k-2, \underbrace{k-1, k-1, \ldots, k-1}_{(a-k+1)\ \textrm{times}}\Bigg), \]
\[ R(0) = \Bigg(1, 2, \ldots, k-2, k-1, \underbrace{a, a, \ldots, a}_{(b-k+1)\ \textrm{times}}\Bigg). \]
%%%
Note that, if at any point $\sigma(R_i) = \sigma(R_j)$, then $U\sigma(R_i) = U\sigma(R_j)$, because $R_i$ and $R_j$ have %the same number of chips and 
the same neighbors. So, $U\sigma_{2k}$ is represented by
\[ L(1) = \Bigg(b-k+1, b-k+2, \ldots, b-2, b-1,  \underbrace{b, b, \ldots, b}_{(a-k+1)\ \textrm{times}}\Bigg), \]
\[ R(1) = \Bigg(1, 2, \ldots, k-2, k-1, \underbrace{0, 0, \ldots, 0}_{(b-k+1)\ \textrm{times}}\Bigg), \]
%%%
and $U^2\sigma_{2k}$ is represented by
\[ L(1) = \Bigg(b-k+1, b-k+2, \ldots, b-2, b-1,  \underbrace{0, 0, \ldots, 0}_{(a-k+1)\ \textrm{times}}\Bigg), \]
\[ R(1) = \Bigg(a-k+2, a-k+3, \ldots, a-1, a, \underbrace{a-k+1, a-k+1, \ldots, a-k+1}_{(b-k+1)\ \textrm{times}}\Bigg). \]
%%%

We can see that for $t=2, 3, \ldots, 2k-1$, the vertex in $G$ that fires (has $b$ chips if it is in $L$, or $a$ chips if it is in $R$) in position $U^t\sigma_{2k}$  is
\begin{equation}
\begin{cases}
R_{k - \frac{t}{2}} & \textrm{for}\ t\ \textrm{even},\\
L_{k - \frac{t-1}{2}} & \textrm{for}\ t\ \textrm{odd}.
\end{cases}
\end{equation}

So, after $2k$ steps, every vertex will have fired once, and by Proposition~\ref{prop:alleq}, $U^{2t}\sigma_{2k}=\sigma_{2k}$.

It remains to construct initial positions with period 1 or 2. The trivial game with no chips on any vertex has period 1, while the initial position where each vertex in $L$ has $b$ chips, and each vertex in $R$ has 0 chips, can be easily checked to have period 2. Thus, all periods $i, 1 \le i \le a$, and $2i, 1 \le i \le a$, are achievable.
\end{proof}

Combining the results of Corollary \ref{cor:lea} and \ref{prop:existpos} gives our main theorem.

%%%%%%%%%%% THEOREM %%%%%%%%%%%%%
\begin{thm} A nonnegative integer $p$ is a possible period of a position $\sigma$ of the parallel chip-firing game on $K_{a,b}$ if and only if 
\begin{equation} \label{eq:posper} p \in (\{i\ |\ 1 \le i \le \min(a,b)\} \cup \{2i\ |\ 1 \le i \le \min(a,b)\}). \end{equation}
\end{thm}
\begin{proof}
By Corollary~\ref{cor:lea}, no period lengths may lie outside the sets in \eqref{eq:posper}; and in Proposition~\ref{prop:existpos}, we have constructed positions with all such periods.
\end{proof}

%%%%%%%%%%%%% SECTION 4 SECTION 4 SECTION 4 %%%%%%%%%%%%%
\section{Periods of Games on the Complete $c$-Partite Graph}

We again use the vector notation from above to represent the positions of a parallel chip-firing game on the complete $c$-partite graph $G=K_{a_1, a_2, \ldots, a_c}$ formed by joining the anticliques $S_1, S_2, \ldots, S_c$; let the vertices in $S_i$ be $S_{i,1}, S_{i,2}, \ldots, S_{i,a_i}$ for each $1 \le i \le c$. Without loss of generality, we will assume $a_1 \ge a_2 \ge \ldots \ge a_c$. As above, we represent a position on $G$ by the set of vectors 
$$\Bigg\{S_i(t) = \Bigg(U^t\sigma(S_{i,1}), U^t\sigma(S_{i,2}), \ldots, U^t\sigma(S_{i,a_i})\Bigg)\Bigg\}, \; 1 \le i \le c.$$

%We claim that all periods from $1$ to $\min(2(a_1+a_2+\ldots+a_{c-1}), a_1+a_2+\ldots+a_c = V(G))$ are possible on $K_{a_1,a_2,\ldots,a_c}$. We construct such positions as follows:

Below is a representation of a position which has period $(c-j) a_c - k + 1$ for all $0 \le j \le c-1$ and $1 \le k \le a_c$.
 Note that a vertex in $S_b$ fires when it has at least $d_b = \sum_{i=0}^{c}(S_i) - S_b$ chips; here $d_b$ is the degree of any vertex in $S_b$.
 
For our construction, we let
\begin{eqnarray*}
%S_{c-i} & = & \Bigg((c-1)-1, 2(c-1)-1, \ldots, (a_c-k)(c-1)-1, \\
% & & \, (a_c-k+1)(c-1)-1, (a_c-k+2)(c-1)-2, \ldots , (a_c-k+(k-1))(c-1)-(k-1), \\ 
% & & \qquad \underbrace{d_{i} - k}_{a_i - a_c + 1}\Bigg) 
S_h(0) & = & \Bigg((h-1), (c-1)+(h-1), 2(c-1)+(h-1), \ldots, (a_c-k-1)(c-1)+(h-1) \\
 & & \, (a_c-k)(c-1)+(h-1), (a_c-k+1)(c-1)+(h-1)-1, \\
 & & \, \ldots, (a_c-k+(k-2))(c-1)+(h-1)-(k-2), \underbrace{d_h - k-\sum_{z=1}^{h-1}\left(a_h - a_c\right)}_{\textrm{appears}\; a_h - a_c + 1\ \textrm{times}}\Bigg) \\
 & & \textrm{for}\ 1 \le h \le c-j-1, \\
S_{i}(0) & = & \Bigg(c-1, 2(c-1), \ldots, (a_c-k)(c-1), \underbrace{d_i, d_i, \ldots, d_i}_{a_i - (a_c-k) \ \textrm{times}}\Bigg) \qquad \textrm{for}\ c-j \le i \le c.
\end{eqnarray*}

We now show that this position indeed has period $(c-j) a_c - k + 1$. Let $F_{\sigma}(t)$ be the set of all firing vertices in $U^t(\sigma)$. 

It can be checked that

%\begin{align*}
$$F_{\sigma}(0) = \bigcup_{i=c-j}^{c}\bigcup_{m=a_c-k}^{a_i}\{S_{i,m}\}; $$
$$F_{\sigma}(t) = \bigcup_{m=a_c}^{a_{c-j-t}}\{S_{c-j-t,m}\}$$ 
for $1 \le t \le c-j-1$; and
$$F_{\sigma}(k_1(c-j-1)-t_1) = \{S_{t_1+1,c-k_1}\}$$
for $1\le k_1\le k, 0\le t_1\le c-j-2$, encompassing steps $c-j$ through $k(c-j-1)$;
$$F_{\sigma}(k(c-j-1)+1+k_2(c-j)) = \bigcup_{i=c-j}^{c}\{S_{c-i, a_c-(k+k_2)}\}$$
 for $0\le k_2 \le a_c-k+1$; and
 
$$F_{\sigma}(k(c-j-1)+1+k_2(c-j)+t_2) = \{S_{c-j-t_2, a_c-(k+k_2)}\} $$
 for $0\le k_2 \le a_c-k-1, 1 \le t_2 \le c-j-1$. (In fact, each vertex $v\in G$ fires exactly when it contains $\textrm{deg}(v)$ chips.)
 
The latter two categories describe which vertices fire during steps $k(c-j-1)+1$ through $k(c-j-1) + 1 + (a_c-k-1)(c-j) + (c-j-1) = (k+1)(c-j-1)+1 + (a_c-k-1)(c-j) = a_c(c-j) - k$. But after this $(a_c(c-j) - k) ^{\textrm{th}}$ step, every vertex in $G$ has fired exactly once; the last to fire is $S_{1,1}$. Hence, $U^{a_c(c-j)-k+1}\sigma = U^0\sigma$ by Proposition~\ref{prop:alleq}, and the period is $a_c(c-j)-k+1$ as desired. This means all periods from $1$ to $c \min(a_1, a_2, \ldots, a_c)$ are achievable, as $j$ ranges from $0$ to $c-1$ and $k$ ranges from $1$ to $a_c$.
As an example, consider the following position on the graph $K_{6, 5, 5, 4}$ with period $11$:
\begin{align*}
S_1(0) & = (0, 3, 6, 7, 7, 7) & & d_1 = 14 \\
S_2(0) & = (1, 4, 7, 10, 10) & & d_2 = 15 \\
S_3(0) & = (3, 6, 15, 15, 15) & & d_3 = 15 \\
S_4(0) & = (3, 6, 16, 16) & & d_4 = 16
\end{align*}

For this position, $a_c=4, c=4, j=1, k=2$, and its predicted period length is $(c-j)a_c - k + 1=11$ as desired.

\section{Discussion and Further Work}

For several graphs, a proof of Bitar's conjecture that $p(\sigma) \le |V(G)|$ for all parallel chip-firing games on those graphs would be interesting; we proved the conjecture for the complete bipartite graph. Though we have constructed many periods of games on complete $c$-partite graphs in Section 4, there exist periods longer than those detailed. For example, take the following position on $K_{2,2,1}$, which has period 5:
\begin{align*}
S_1(0) & = (2) & & d_1 = 4 \\
S_2(0) & = (1, 2) & & d_2 = 3 \\
S_3(0) & = (0, 3) & & d_3 = 3
\end{align*}
Though positions with these larger periods are more difficult to characterize generally, Bitar's conjecture still appears to be true for complete $c$-partite graphs. 

 Moreover, bounding the periods of positions on vertex-regular graphs and more general bipartite graphs are directions for further research. By doubling the length of each cycle in the graph used in the counterexample by Kiwi et.\ al.\ \cite{kiwi}, we find a counterexample on a graph containing only even cycles, that is, for the general bipartite graph. 
 
We would also like to determine which periods less than the bound are possible. Levine \cite{levine} related period lengths of games on the complete graph to the \emph{activity}, defined as $\displaystyle \lim_{t\rightarrow \infty} \frac{\sum_{v \in G} u_t(\sigma, v)}{vt}$. On the other hand, we believe that period lengths are related to lengths of subcycles (closed paths) of the graph $G$; in particular, we conjecture that any period length of a game on $G$ is either a divisor of the order of some subcycle of $G$, or perhaps the least common multiple of the orders of some disjoint subcycles of $G$. This agrees with known results for the tree graph \cite{bitar}, complete graph \cite{levine}, and now the complete bipartite graph. Our numerical experiments have also verified this conjecture for cycle graphs and complete $k$-partite graphs; in fact, my correspondence with Zhai \cite{zhai} has produced a proof of this conjecture for the cycle graph.

Another interesting direction to pursue is observing the implications of ``reducing'' the parallel chip-firing game by removing as many chips as possible from each vertex without affecting their firing pattern (without changing $F_v(t)$ for all $v \in G$ and $t\ge 0$). This reduction may simplify some games into being more approachable by induction.

Besides studying period lengths of parallel chip-firing games, an examination of the transient length of games on certain graphs would be useful in modeling real-world phenomena. Studying transient positions would also help uncover what attributes determine whether a position is within a period or not, and bounding the transient length would make for more efficient computation of the period length of games on complex graphs.

Chip-firing games on lattices and tori have been used as cellular automaton models of the deterministic fixed-energy sandpile (see \cite{casartelli,bagnoli}). Since most studies of sandpiles have been concerned with asymptotic measures such as the ``activity,'' bounding the period length of such models could serve as a measure of the fidelity of the model to the real world.
%For acknowledgements section, please don't number the section, please begin it with \section*{Acknowledgements}

\section*{Acknowledgments} 
The author would like to thank my mentor, Yan Zhang of the Massachusetts Institute of Technology, for his teaching, guidance, and support, the Center for Excellence in Education and the Research Science Institute for sponsoring my research, and Dr.\ John Rickert for tips on writing and lecturing. The author would also like to thank Dr.\ Lionel Levine, Ziv Scully, Scott Kominers, Paul Kominers, Daniel Vitek, Brian Hamrick, Dr.\ Dan Teague, Patrick Tenorio, Alex Zhai, and Dr.\ Ming Ya Jiang for valuable pointers and advice on this problem and this paper.

% You may incorporate your references as follows in your main tex file.
% Using BibTex is not recommended but can be handled.

\medskip
% The data information below will be filled by AIMS editorial staff
\medskip

\end{document}